\documentclass[12pt]{amsart}
\usepackage{amscd}
\def\NZQ{\Bbb}               
\def\NN{{\NZQ N}}

\def\FFF{{\NZQ F}}
\def\LLL{{\NZQ L}}
\def\JJ{{\mathcal J}}
\def\LL{{\mathcal L}}
\def\frk{\frak}               

\def\mm{{\frk m}}

\def\Phi{{\frk n}}
\def\Phi{{\frk N}}
\def\opn#1#2{\def#1{\operatorname{#2}}} 
\opn\chara{char}
\opn\length{\ell}
\opn\pd{pd}
\opn\rk{rk}
\opn\projdim{proj\,dim}
\opn\injdim{inj\,dim}
\opn\rank{rank}
\opn\depth{depth}
\opn\grade{grade}
\opn\height{height}
\opn\embdim{emb\,dim}
\opn\codim{codim}

\opn\Tr{Tr}
\opn\bigrank{big\,rank}
\opn\superheight{superheight}\opn\lcm{lcm}
\opn\trdeg{tr\,deg}%
\opn\reg{reg}
\opn\lreg{lreg}
\opn\ini{in}
\opn\lpd{lpd}
\opn\div{div}
\opn\Div{Div}
\opn\cl{cl}
\opn\Cl{Cl}
\opn\Spec{Spec}
\opn\Supp{Supp}
\opn\supp{supp}
\opn\Sing{Sing}
\opn\Ass{Ass}
\opn\Ann{Ann}
\opn\Rad{Rad}
\opn\Soc{Soc}
\opn\Im{Im}
\opn\Ker{Ker}
\opn\Coker{Coker}
\opn\Am{Am}
\opn\Hom{Hom}
\opn\Tor{Tor}
\opn\Ext{Ext}
\opn\End{End}
\opn\Aut{Aut}
\opn\id{id}

\opn\nat{nat}
\opn\pff{pf}
\opn\Pf{Pf}
\opn\GL{GL}
\opn\SL{SL}
\opn\mod{mod}
\opn\ord{ord}
\opn\Gin{Gin}
\opn\Hilb{Hilb}
\opn\aff{aff}
\opn\con{conv}
\opn\relint{relint}
\opn\st{st}
\opn\lk{lk}
\opn\cn{cn}
\opn\core{core}
\opn\vol{vol}
\opn\dist{dist}
\opn\link{link}
\opn\star{star}
\opn\gr{gr}
\def\Rees{{\mathcal R}}
\def\pot#1#2{#1[\kern-0.28ex[#2]\kern-0.28ex]}

\opn\dirlim{\underrightarrow{\lim}}
\opn\inivlim{\underleftarrow{\lim}}
\let\union=\cup

\let\Union=\bigcup

\let\Dirsum=\bigoplus

\let\to=\rightarrow
\let\To=\longrightarrow
\def\Implies{\ifmmode\Longrightarrow \else
        \unskip${}\Longrightarrow{}$\ignorespaces\fi}
\def\implies{\ifmmode\Rightarrow \else
        \unskip${}\Rightarrow{}$\ignorespaces\fi}
\def\iff{\ifmmode\Longleftrightarrow \else
        \unskip${}\Longleftrightarrow{}$\ignorespaces\fi}

\let\:=\colon
\newtheorem{Theorem}{Theorem}[section]
\newtheorem{Lemma}[Theorem]{Lemma}
\newtheorem{Corollary}[Theorem]{Corollary}
\newtheorem{Proposition}[Theorem]{Proposition}

\let\epsilon\varepsilon
\let\phi=\varphi
\let\kappa=\varkappa
\def\qed{\ifhmode\textqed\fi
      \ifmmode\ifinner\quad\qedsymbol\else\dispqed\fi\fi}
\def\textqed{\unskip\nobreak\penalty50
       \hskip2em\hbox{}\nobreak\hfil\qedsymbol
       \parfillskip=0pt \finalhyphendemerits=0}
\def\dispqed{\rlap{\qquad\qedsymbol}}

\opn\dis{dis}
\def\pnt{{\raise0.5mm\hbox{\large\bf.}}}

\opn\Lex{Lex}

\begin{document}

\title{The depth of powers of an  ideal}

\author{J\"urgen Herzog and Takayuki Hibi}
\thanks{This paper was completed while the authors stayed at the Mathematisches Forschungsinstitut in Oberwolfach in the frame of the Research in Pairs Program}
\address{J\"urgen Herzog, Fachbereich Mathematik und
Informatik, Universit\"at Duisburg-Essen, Campus Essen,
45117 Essen, Germany}
\email{juergen.herzog@uni-essen.de}

\address{Takayuki Hibi, Department of Pure and Applied Mathematics,
Graduate School of Information Science and Technology,
Osaka University, Toyonaka, Osaka 560-0043, Japan}
\email{hibi@math.sci.osaka-u.ac.jp}

\maketitle

\begin{abstract}
We study the limit  and initial behavior of the numerical function $f(k)=\depth S/I^k$. General properties of this function together with concrete examples arising  from combinatorics are discussed.
\end{abstract} 

\section*{Introduction}

Let $S$ be either a Noetherian local ring with maximal ideal $\mm$, or a standard graded $K$-algebra with graded maximal ideal $\mm$, where $K$ is any field, and let $I\subset S$ be a proper ideal, which we assume to be graded if $S$ is standard graded. We are interested in the behavior of the numerical function $\depth S/I^k$. It is clear that this function is bounded by the dimension $d$ of $S$. A classical result by Burch \cite{Burch} says that 
\[
\min_k\depth S/I^k\leq d-\ell(I),
\]
where $\ell(I)$ the analytic spread of $I$, that is, the dimension of  $\mathcal{R}(I)/\mm\mathcal{R}(I)$. Here  $ \mathcal{R}(I)=\Dirsum_kI^kt^k$ is the Rees ring of $I$.

By a  theorem of Brodmann \cite{Brodmann2}, $\depth S/I^k$ is constant for $k\gg 0$. We call this constant value the {\em limit depth of $I$}, and denote it by $\lim_{k\to\infty}\depth S/I^k$. Brodmann  improved the Burch inequality by showing that 
\[
\lim_{k\to\infty}\depth S/I^k\leq d-\ell(I),
\]
Eisenbud and Huneke \cite{EisenbudHuneke} showed that equality holds, if the associated graded ring $\gr_I(S)$ is Cohen--Macaulay. This is for example the case if $S$ and $\mathcal{R}(I)$ are Cohen--Macaulay, see Huneke \cite{Huneke}.
Recently  Branco Correia and Zarzuela \cite{CorreiaZarzuela} proved similar results for Rees powers of a module. 
In Section 1 we will  give new and relatively short proofs for these facts. 

While the limit behavior of $\depth S/I^k$ is well understood, the initial behavior of $S/I^k$ is more mysterious. On the  one  hand,  if one chooses a homogenous ideal `randomly', one can be quite sure that $\depth S/I^k$ is a decreasing function. So this behavior seems to be the normal one. On the other hand,  Trung and Goto independently communicated to the first author examples of graded ideals such that $S/I^2$ is Cohen--Macaulay, while $S/I$ is not Cohen--Macaulay. In these cases, of course, $\depth S/I<\depth S/I^2$.

In Section 2 we show that $\depth S/I^k$ is a decreasing function if all powers of $I$ have a linear resolution, and we show that all powers of a monomial ideal have linear quotients, and hence have linear resolutions, if with respect to a suitable monomial order, the toric  ideal $J$ of the Rees ring of $I$  satisfies the so-called $x$-condition, which is a condition on the Gr\"obner basis  of $J$. If this condition is satisfied, one also obtains lower bounds for  $\depth S/I^k$. We also derive a formula for $\depth S/I$ when $I$ has linear quotients. 

We use the techniques developed in the first sections to compute the function $\depth S/I^k$ for classes of ideals arising in combinatorial contexts. By \cite{HHZpower} we know that the $x$-condition is satisfied for all edge ideals of finite graphs whose complementary graph is chordal. Thus all powers of such ideals have linear quotients. 

We next consider polymatroidal ideals. Powers of polymatroidal ideals are again polymatroidal. Since polymatroidsal ideals have linear quotients we can compute depth $S/I^k$ for all $k$. Explicit formulas are given for special classes of polymatroidal ideals, namely for ideals of Veronese type. 

Finally we consider monomial  ideals coming from  finite posets. In this case, again all powers have linear quotients. Choosing posets suitably we can show that, given a decreasing function $f : \NN \to \NN$ with 
$f(0) = 2 \lim_{k \to \infty} f(k) + 1$
for which $\Delta f$ is decreasing, there exists a monomial ideal $I
\subset S$ such that $\depth S / I^k = f(k)$ for all $k \geq 1$. Here 
$(\Delta f)(k) = f(k) -
f(k+1)$ for all $k \in \NN$. 

All examples considered in Section 3 have decreasing depth functions. However we show in Section 4 that, given any bounded increasing numerical function $f : \NN \setminus \{ 0 \} \to \NN$, there exists a monomial ideal $I$ such that $\depth S/I^k=f(k)$ for all $k$. In all cases mentioned so far, the depth function is monotonic. We conclude this paper with an example of a monomial ideal whose depth function is not monotonic. 

In  view of the examples in this paper we are tempted to conjecture that the depth function can be any convergent numerical nonnegative function.

\section{The limit  behavior of $\depth S/I^k$}

Let  $A$ be finitely generated a standard graded $R$-algebra, and  $E$ be a finitely generated graded $A$-module. Then each graded component $E_k$ of $E$ is a finitely generated $R$-module. 

\begin{Theorem}
\label{basic}
The depth of $E_k$ is constant for $k\gg 0$, and hence $\lim_{k\to \infty}\depth E_k$ exists. Moreover one has  $$\lim_{k\to \infty}\depth E_k\leq \dim E-\dim E/\mm E,$$
and equality holds if $E$ is Cohen--Macaulay. 
\end{Theorem}

\begin{proof}
Let $x_1,\ldots,x_n$ be a minimal set of generators of $\mm$. Then $\depth E_k=n-\max\{i\: H_i(x;E_k)\neq 0\}$, see \cite{BH}. Here we denote by $H(x;M)$ the Koszul homology of a module $M$ with respect to the sequence $x=x_1,\ldots, x_n$.

Now consider the homology modules $H_i(x;E)$. These modules a finitely generated graded $A$-modules with graded components 
\[
H_i(x;E)_k=H_i(x;E_k).
\]
Let $c=\max\{i\: \dim H_i(x;E)>0\}$. Then for all $i>c$, we have $\dim H_i(x;E)=0$, so that $H_i(x;E)_k=0$ for all $i>c$ and all $k\gg 0$. On the other hand, since $\dim H_c(x;E)>0$ it follows that $H_c(x;E)_k\neq 0$ for all $k\gg 0$.  This implies that $\depth E_k=n-c$ for all $k\gg 0$. 

Let $E^{(r)}=\Dirsum_iE_{ir}$ the  $r$th-Veronese submodule of $E$. Note that $\dim E^{(r)}=\dim E$, $\dim E/\mm E=\dim E^{(r)}/\mm E^{(r)}$, and that $\depth E^{(r)}_k=\depth E_{rk}$ is constant for $r\gg 0$ and all $k$. Thus if we replace $E$ by $E^{(r)}$ for $r$ big enough, we may assume that 
\[
\grade(\mm, E)=n-\max\{i\: H_i(x;E)\neq 0\}=\lim_{k\to\infty}\depth E_k.
\]
Since $\grade(\mm,E)\leq \dim E-\dim E/\mm E$ with equality if $E$ is Cohen--Macaulay (see \cite[Theorem 2.1.2]{BH}), the assertions follow.  
\end{proof}

As a consequence we obtain the  theorem of Brodmann \cite{Brodmann2} together with a statement on $\lim_{k\to\infty} \depth I^k$, as well as the result of Eisenbud and Huneke \cite{EisenbudHuneke}.

\begin{Theorem}
\label{powers}
The limits $\lim_{k\to \infty}\depth I^k$, $\lim_{k\to \infty} \depth S/I^k$ and\newline  $\lim_{k\to \infty}\depth  I^k/I^{k+1}$ exist, and 
\[
\lim_{k\to \infty}\depth S/I^k\leq \lim_{k\to\infty} \depth I^k-1=\lim_{k\to\infty}\depth I^k/I^{k+1}\leq  \dim S-\ell(I).
\]
If $S$ is Cohen--Macaulay and $\height I>0$, then $$\lim_{k\to \infty}\depth S/I^k= \lim_{k\to\infty} \depth I^k-1.$$ Moreover,   all limits are equal to $\dim S-\ell(I)$ if, in addition, the associated graded ring $\gr_I(S)$ is Cohen-Macaulay.
\end{Theorem}

\begin{proof} Let $E=\mathcal{R}(I)$ the Rees ring, or $E=\gr_I(S)$ the associated graded ring of $I$. In the first case, Theorem \ref{basic} implies that $\lim_{k\to\infty}\depth  I^k$ exists, in the second case the theorem  implies that 
$\lim_{k\to\infty}\depth I^k/I^{k+1}$ exists.

The last inequality also follows  from Theorem \ref{basic}, since $\dim \gr_I(S)=\dim S$ and $\dim  \gr_I(R)/\mm\gr_I(S)= \dim \mathcal{R}(I)/\mm\mathcal{R}(I)=\ell(I)$.

Now we show that $\lim_{k\to \infty}\depth S/I^k$ exists. To this end we consider the exact sequences
\[
0\To I^k/I^{k+1}\To S/I^{k+1}\To S/I^k\To 0.
\]
Set $f(k)=\depth S/I^k$, and let $c=\lim_{k\to \infty}I^k/I^{k+1}$. Then there exists an integer $k_0$ such that  for all $k\geq k_0$ these exact sequences give rise to the following inequalities
\begin{enumerate}
\item[(i)] $f(k+1)\geq \min\{c,f(k)\}$,
\item[(ii)] $c\geq \min\{f(k+1), f(k)+1\}$,
\end{enumerate}
see \cite[Proposition 1.2.9]{BH}. Suppose  that $f(k)\geq  c$ for some  $k\geq k_0$. Then (ii) implies that $f(k+1)\leq c$. Then (i) yields that $f(k+1)=c$. It follows that $f(\ell)=c$ for all $\ell\geq k+1$. Hence $\lim_{k\to\infty}f(k)=c$ in this case. 

We may henceforth assume that $f(k)\leq c$ for all $k$. 
Then (i) implies that $f(k)$ is an increasing function for $k\geq k_0$, and that this function  is bounded above by $c$. Thus the limit $f(k)$ exists, and it is less than or equal to $c$.

Next we want to prove the equation $\lim_{k\to\infty}\depth I^k-1=\lim_{k\to\infty}\depth I^k/I^{k+1}$. 
The short exact sequence 
\[
0\to I^k\to I^{k+1}\to I^k/I^{k+1}\to 0
\]
yields for $k\geq k_0$ the inequalities
\[
c\geq \min\{g(k+1)-1,g(k)\},
\]
where $g(k)=\depth I^k$. Let $g=\lim_{k\to \infty} g(k)$. Then passing to the limit we see that $c\geq \min\{g-1,g\}=g-1$. 

Suppose $c>g-1$, and let  $n$ be the minimal number of generators of $\mm$. Then there exists an integer $k_0$ such that  $H_{n-g}(x;I^k)\neq 0$ and $H_{n-g+1}(x;I^k/I^{k+1})=0$ for all $k\geq k_0$. This implies that  the natural map $H_{n-g}(x;I^{k+1})\to H_{n-g}(x;I^k)$ is injective for all $k\geq k_0$. Composing these maps, we see that $H_{n-g}(x;I^{\ell})\to H_{n-g}(x;I^k)$ are injective for $k\geq k_0$ and all $\ell>k$. However, the Artin--Rees lemma implies that for any finitely generated $S$-module, the natural homomorphism $H_{n-g}(x;I^{\ell}M)\to H_{n-g}(x;M)$ is the zero map for $\ell\gg 0$. Thus we conclude that $H_{n-g}(x;I^{\ell})=0$ for $\ell\gg 0$, a contradiction.

Suppose now that $S$ is Cohen-Macaulay, and that $\height I>0$. Then $\depth S/I^k=\depth I^k-1$, so that 
$\lim_{k\to\infty}S/I^k=\lim_{k\to\infty}I^k-1$. Finally, if $\gr_I(S)$ is Cohen-Macaulay, then $\lim_{k\to\infty}I^k/I^{k+1}=\dim S-\ell(I)$, by Theorem \ref{basic}. \hspace{2cm}
\end{proof}

\section{The initial behavior of $\depth S/I^k$}

On support of the normal behavior we show 

\begin{Proposition}
\label{linear}
 Let $I$ be a graded ideal all of whose  powers have a linear resolution. Then $\depth S/I^k$ is a decreasing function of $k$.
\end{Proposition}

The proposition is a consequence of Corollary \ref{increasing} stated below. As usual we denote by $\beta_{ij}(M)$ the graded Betti numbers of a graded module $M$ over $S$. We call the least degree of homogeneous generator of $M$, the {\em initial degree of $M$.}

\begin{Lemma}
\label{nothing}
Let $J\subset I$ be graded ideals, and let $d$ be the initial degree of $I$. Then 
\[
\beta_{i,i+d}(J)\leq \beta_{i,i+d}(I)
\]
for all $i$.
\end{Lemma}

\begin{proof}
The short exact sequence 
\[
0\To J\To I\To I/J\To 0
\] 
yields the long exact sequence
\[
\cdots \To \Tor_{i+1}(K,I/J)_{i+1+(d-1)}\To \Tor_{i}(K,J)_{i+d}\To\Tor_{i}(K,I)_{i+d}\To\cdots
\]

Since  the initial degree of $I/J$ is greater than or equal to $d$, it follows that 
$\Tor_{i+1}(K,I/J)_{i+1+(d-1)}=0$. Hence $\Tor_{i}(K,J)_{i+d}\to \Tor_{i}(K,I)_{i+d}$ is injective.
\end{proof}

Let $\FFF$ be the graded  minimal free resolution of $I$, and suppose that $d$ is the initial degree of $I$. Then the subcomplex $\LLL$ of $\FFF$ with $L_i=S(-i-d)^{\beta_{i,i+d}}$ is called the {\em lowest linear strand of $\FFF$}. We call its length the {\em linear projective dimension of $I$}.

\begin{Corollary} 
\label{increasing}
Let $I\subset S$ be a graded ideal with initial degree $d$. Then $$\beta_{i,i+(k+1)d}(I^{k+1})\geq \beta_{i,i+kd}(I^k)$$ for all $k$. In particular, the linear projective dimension of $I^k$ is an increasing function of $k$. 
\end{Corollary}

\begin{proof}
Let $x\in I\setminus\mm I$. Then $xI^k\subset I^{k+1}$. It follows from Lemma \ref{nothing} that $\beta_{i,i+kd}(I^k)=\beta_{i,i+(k+1)d}(xI^k)\leq \beta_{i,i+{(k+1)}d}(I^{k+1})$.
\end{proof}

We now discuss graded ideals having linear quotients.
Let $f_1, \ldots, f_s$ be a sequence of homogeneous elements
of $S$ with $0<\deg f_1\leq \deg f_2\leq\cdots\leq \deg f_s$. 
We say that $f_1, \ldots, f_s$
has {\em linear quotients} if, 
for each $2 \leq j \leq s$, 
the colon ideal $(f_{1}, f_{2}, \ldots, f_{j-1}):f_{j}$
is generated by linear forms.
We say that a graded ideal $I \subset S$ has linear quotients
if $I$ is generated by a sequence with linear quotients.
It is known \cite[Lemma 4.1]{ConcaHerzog} that 
if $f_1, \ldots, f_s$ is a sequence with linear quotients
and if all $f_i$ have the same degree, then the ideal
$(f_1, \ldots, f_s)$ has a linear resolution.  

Let $I$ be a graded ideal generated by a sequence
$f_1, \ldots, f_s$ with linear quotients.
Let $q_j(I)$ denote the minimal number of linear forms
generating $(f_{1}, f_{2}, \ldots, f_{j-1}):f_{j}$,
and $q(I) = \max_{2 \leq j \leq s}q_j(I)$.

As in the proof of
\cite[Corollary 1.6]{HerzogTakayama}
we can show that 
the length of the minimal free resolution of 
$S/I$ over $S$ is equal to $q(I) + 1$.
Hence  
\begin{eqnarray}
\label{depth}
\depth S/I = n - q(I) - 1.
\end{eqnarray}
Thus in particular the integer $q(I)$ is independent of 
the particular choice of the sequence 
of generators with linear quotients.

\begin{Corollary}
\label{moregeneral}
Let $I$ be a graded ideal generated in degree $d$ with linear resolution, and let $f_1,\ldots, f_s$ be a sequence with linear quotients which is part of a minimal system of generators of $I$. Then $\depth S/I\geq n-q(J)-1$, where $J$ is the ideal generated by $f_1,\ldots, f_s$.
\end{Corollary}

\begin{proof}
Since $J$ and $I$ both have a linear resolution, it follows from Lemma \ref{nothing} that $\projdim S/J\leq \projdim S/I$. Hence $\depth S/I\geq \depth S/J=n-q(J)-1$.
\end{proof}

\medskip
\noindent
Our next goal is to discuss a  Gr\"obner basis condition that guarantees that all powers of an ideal have linear quotients. Let $K$ be a field and $S = K[x_1, \ldots, x_n]$ the polynomial ring
in $n$ variables over $K$ with each $\deg x_i = 1$.  Let $I \subset
S$ be a monomial ideal generated in one degree and $G(I)$ its minimal
system of monomial generators.  Recall that the Rees algebra
$\Rees(I)$ of $I$ is 
\[
\Rees(I) = K[x_1, \ldots, x_n, \{ ut \}_{u \in G(I)}] \subset S[t].
\]
Let $A = K[x_1, \ldots, x_n, \{ y_u \}_{u \in G(I)}]$ denote the
polynomial ring in $n + |G(I)|$ variables over $K$ with each $\deg
x_i = \deg y_u = 1$.  The {\em toric ideal} of $\Rees(I)$ is the
kernel $J_{\Rees(I)}$ of the surjective homomorphism $\pi : A \to
\Rees(I)$ defined by setting $\pi(x_i) = x_i$ for all $1 \leq i \leq
n$ and $\pi(y_u) = ut$ for all $u \in G(I)$.

Let $<_{lex}$ denote the lexicographic order on $S$ induced by $x_1 >
x_2 > \cdots > x_n$.  Fix an arbitrary monomial order $<^{\#}$ on
$K[\{ y_u \}_{u \in G(I)}]$.  We then introduce the new monomial
order $<^{\#}_{lex}$ on $A$ defined as follows:
For monomials
$(\prod_{i=1}^{n} x_i^{a_i})(\prod_{u \in G(I)} y_{u}^{a_u})$ and 
$(\prod_{i=1}^{n} x_i^{b_i})(\prod_{u \in G(I)} y_{u}^{b_u})$
belonging to $A$, one has
\[
(\prod_{i=1}^{n} x_i^{a_i})(\prod_{u \in G(I)} y_{u}^{a_u})
<^{\#}_{lex}
(\prod_{i=1}^{n} x_i^{b_i})(\prod_{u \in G(I)} y_{u}^{b_u})
\]
if either
\begin{enumerate}
\item[(i)]
$\prod_{u \in G(I)} y_{u}^{a_u} <^{\#} \prod_{u \in G(I)} y_{u}^{b_u}$
or 
\item[(ii)] 
$\prod_{u \in G(I)} y_{u}^{a_u} = \prod_{u \in G(I)} y_{u}^{b_u}$
and $\prod_{i=1}^{n} x_i^{a_i} <_{lex} \prod_{i=1}^{n} x_i^{b_i}$.
\end{enumerate}
Let ${\mathcal G}(J_{\Rees(I)})$ denote the reduced Gr\"obner
basis of $J_{\Rees(I)}$ with respect to $<^{\#}_{lex}$. We say that $I$ satisfies the {\em $x$-condition} if 
each element belonging to
${\mathcal G}(J_{\Rees(I)})$ is at most linear in the variables
$x_1, \ldots, x_n$.

\begin{Theorem}
\label{atmostlinear}
Suppose that $I$ satisfies the $x$-condition.  Then each power of $I$ has linear quotients.
\end{Theorem}

\begin{proof}
Fix $k \geq 1$.  Each $w \in G(I^k)$ has a unique expression, 
called the {\em standard expression}, of the form 
$w = u_1 \cdots u_k$ 
with each $u_i \in G(I)$ such that 
$y_{u_1} \cdots y_{u_k}$ is a standard monomial of $A$ 
with respect to $<^{\#}$, that is, a monomial which does not belong to the initial ideal of $J_{\Rees(I)}$.  Let $w^{*}$ denote
the standard monomial $y_{u_1} \cdots y_{u_k}$.  
Let $G(I^k) = \{ w_1, \ldots, w_s \}$ with 
$w_1^{*} <^{\#} \cdots <^{\#} w_s^{*}$.

We claim that $I^k$ has linear quotients with the ordering 
$w_1, \ldots, w_s$ of its generators. 
Let $f$ be a monomial belonging to the colon ideal
$(w_1, \ldots, w_{j-1}) : w_j$. 
Thus $fw_j = gw_i$ for some $i < j$ and for some monomial $g$. 
Let $w_j = u_1 \cdots u_k$ and $w_i = v_1 \cdots v_k$ be the standard
expressions of $w_j$ and $w_i$. 
The binomial 
$f y_{u_1} \cdots y_{u_s} - g y_{v_1} \cdots y_{v_s}$
belongs to $J_{\Rees(I)}$. 
Since
$y_{v_1} \cdots y_{v_s} <^{\#} y_{u_1} \cdots y_{u_s}$,
it follows that the initial monomial of $f y_{u_1} \cdots y_{u_s} - g y_{v_1} \cdots y_{v_s}$ is 
$f y_{u_1} \cdots y_{u_s}$.
Hence there is a binomial 
$h^{(+)} - h^{(-)}$
belonging to ${\mathcal G}(J_{\Rees(I)})$ 
whose initial monomial $h^{(+)}$ divides
$f y_{u_1} \cdots y_{u_s}$.
Since $y_{u_1} \cdots y_{u_s}$ is a standard monomial 
with respect to $<^{\#}$, it follows from
the definition of the monomial order $<_{lex}^{\#}$
that it remains to be a standard monomial 
with respect to $<_{lex}^{\#}$.  
Hence the initial monomial of none of the binomials belonging to 
${\mathcal G}(J_{\Rees(I)})$ can divide  
$y_{u_1} \cdots y_{u_s}$.  As a consequence, 
the initial monomial $h^{(+)}$ 
must be divided by some variable, say, $x_a$.
Since $h^{(+)}$ is at most linear in the variables
$x_1, \ldots, x_n$, one has
$h^{(+)} = x_a y_{u_{p_1}} \cdots y_{u_{p_t}}$; then
$x_a$ divides $f$ and where
$y_{u_{p_1}} \cdots y_{u_{p_t}}$ divides
$y_{u_1} \cdots y_{u_s}$.
Let $h^{(-)} = x_b y_{v_{q_1}} \cdots y_{v_{q_t}}$, where 
$y_{v_{q_1}} \cdots y_{v_{q_t}} <^{\#} y_{u_{p_1}} \cdots
y_{u_{p_t}}$.
One has
$x_a u_{p_1} \cdots u_{p_t} = x_b v_{q_1} \cdots v_{q_t}$.

To complete our proof, 
we show that  $x_a \in (w_1, \ldots, w_{j-1}) : w_j$.
Since $y_{u_{p_1}} \cdots y_{u_{p_t}}$ divides
$y_{u_1} \cdots y_{u_s}$, we can write
$y_{u_1} \cdots y_{u_s}
= y_{u_{p_1}} \cdots y_{u_{p_t}} y_{u_{p_{t+1}}} \cdots y_{u_{p_k}}$.

Since 
$y_{v_{q_1}} \cdots y_{v_{q_t}} <^{\#}
y_{u_{p_1}} \cdots y_{u_{p_t}}$,
it follows that
$$y_{v_{q_1}} \cdots y_{v_{q_t}} y_{u_{p_{t+1}}} \cdots y_{u_{p_k}}
<^{\#} y_{u_1} \cdots y_{u_k} = w_j^{*}.$$
Let 
$w_{i_0} = v_{q_1} \cdots v_{q_t} u_{p_{t+1}} \cdots u_{p_k} \in
G(I^k)$.
Then $x_a w_j = x_b w_{i_0}$.
Since $w_{i_0}^* \leq^{\#} 
y_{v_{q_1}} \cdots y_{v_{q_t}} y_{u_{p_{t+1}}} \cdots y_{u_{p_k}}$, 
one has $w_{i_0}^* < ^{\#} w_j^{*}$.  Hence $i_0 < j$.  Thus $x_a
\in (w_1, \ldots, w_{j-1}) : w_j$, as desired.
\end{proof}

 We write $\ini(J_{\Rees(I)})$ for the initial ideal of $J_{\Rees(I)}$ with respect to the monomial order $<_{lex}^{\#}$  introduced above.

Let $m=|G(I)|$. For each multi-index $a=(a_1,\ldots, a_m)\in \NN^m$, we set $|a|=\sum_{i=1}^ma_i$.

\begin{Corollary}
\label{schappach}
Let 
\[
\rho(a)=|\{i\: x_iy^a\in \ini(J_{\Rees(I)})\}|,
\]
and let $c$ be the least common multiple of those monomials $y^a$ for which there exits an integer $i$ such that $x_iy^a\in G(\ini(J_{\Rees(I)}))$. Then 
\begin{enumerate}
\item[(a)] $\depth S/I^k\geq n-\max\{\rho(a)\: |a|=k\}-1$;
\item[(b)] $\lim_{k\to\infty}\depth S/I^k\geq n-\rho(c)-1$.
\end{enumerate}
\end{Corollary}

\begin{proof}
We consider $A=S[y_1,\ldots, y_m]$ a bigraded $K$-algebra with each $\deg x_i=(1,0)$ and each $\deg y_j=(0,1)$. 
Then $J=J_{\Rees(I)}$ is a  bigraded ideal. For each $k$,  $J_{(*,k)}=\Dirsum_{i}J_{(i,k)}$ is a submodule of the free $S$-module  $A_{(*,k)}=\Dirsum_{a, \; |a|=k}Sy^a$, and one has the free presentation 
\[
0\to J_{(*,k)}\To A_{(*,k)}\To I^k\To 0.
\]
On the free $S$-module $A_{(*,k)}$ we introduce the monomial order induced by the monomial order $<_{lex}^{\#}$. Then we have 
\[
\ini(J_{(*,k)})=\ini(J)_{(*,k)}.
\]
By a standard deformation argument it follows  therefore 
\[
\projdim I^k\leq \projdim A_{(*,k)}/\ini(J)_{(*,k)}.
\]
We have 
\[
\ini(J)_{(*,k)}=\Dirsum_{a,\; |a|=k}L_ay^a,
\]
where $L_a$ is generated by all $x_i$ such that $x_iy^a\in \ini(J)$. Therefore $$\projdim \ini(J)_{(*,k)}= 
\max\{ \rho(a)\: |a|=k\}-1.$$ Thus assertion (a) follows. 
Statement (b) is a simple consequence of (a), observing that $L_a\subset L_b$ if $y^a$ divides $y^b$.
\end{proof}

\section{Classes of examples arising in combinatorics}

The function $\depth S/I^k$ will be computed for certain classes of monomial ideals, viz., polymatroidal ideals,
edge ideals of finite graphs, and monomial ideals of finite lattices.

\medskip
\noindent
(a) A typical example for which Theorem \ref{atmostlinear} can be applied
arises from a finite graph.
Let $G$ be a finite graph on the vertex set 
$[n] = \{ 1, \ldots, n \}$,
having no loop and no multiple edge, with $E(G)$ its edge set.  
Let, as before, 
$S = K[x_1, \ldots, x_n]$ denote the polynomial ring 
in $n$ variables over $K$. 
The {\em edge ideal} of $G$ is the ideal $I(G)$ of $S$ which is
generated by those 
quadratic monomials $x_i x_j$ with $\{ i, j \} \in E(G)$.  
It is known \cite{Froeberg} that $I(G)$ has a linear resolution 
if and only if 
the complementary graph of $G$ is chordal.  
(Recall that the {\em complementary graph} of $G$ is the finite 
graph ${\bar G}$ with 
$E({\bar G}) 
= \{ \{ i, j \} \subset [n] : \{ i, j \} \not\in E(G) \}$.  
On the other hand, a finite graph is called {\em chordal} 
if each of its cycles of length $\geq 4$ has a chord.)  
Moreover, in \cite[Theorem 3.2]{HHZpower} it is proved that 
if ${\bar G}$ is chordal, then $I(G)$ has linear quotients.

In the following we assume that $\bar{G}$ is chordal. In this case the clique complex of $\bar{G}$ is a quasi-forest and we order the vertices according to a leaf order of this quasi-forest, see the proof of \cite[Proposition 2.3]{HHZpower}.
Let
$
\Rees(I(G)) = K[x_1, \ldots, x_n, \{ x_i x_j t \}_{\{ i,j \} \in
E(G)}]
$
denote the Rees algebra of $I(G)$,
$A = K[x_1, \ldots, x_n \{ y_{i,j} \}_{\{ i,j \} \in E(G)}$
the polynomial ring in $n + |E(G)|$ variables over $K$, 
and $J_{\Rees(I(G))}$ the toric ideal of $\Rees(I(G))$.  
Thus $J_{\Rees(I(G))}$ is the kernel of the surjective homomorphism 
$\pi : A \to \Rees(I(G))$ defined by setting $\pi(x_i) = x_i$ 
for all $i$ 
and $\pi(y_{ij})= x_ix_j$ for all $\{ i, j \} \in E(G)$. 
We introduce the ordering $<$ of the variables of $A$ by setting
(i) $y_{i,j} > y_{p,q}$, where $i < j$ and $p < q$, if 
either $i < p$ or ($i = p$ and $j < q$), and 
(ii) $y_{i,j} > x_1 > \cdots > x_n$ for all 
$\{ i, j \} \in E(G)$.  
Let $<_{lex}$ denote a lexicographic order on $A$ 
induced by the ordering $<$ and ${\mathcal G}(J_{\Rees(I(G))})$ 
the reduced Gr\"obner basis of $\Rees(I(G))$ 
with respect to $<_{lex}$.

\medskip
We quote the following result \cite[Theorem 3.1]{HHZpower} 

\begin{Theorem}
\label{finitegraph}
Suppose that the complementary graph of $G$ is chordal. Then each element belonging to
${\mathcal G}(J_{\Rees(I(G))})$ is at most linear in the variables
$x_1, \ldots, x_n$.  
\end{Theorem}

In \cite[Theorem 3.2]{HHZpower} it is proved that if 
${\bar G}$ is chordal, then
each power of $I(G)$ has a linear resolution.  
By virtue of Theorem \ref{atmostlinear}, 
we  have

\begin{Corollary}
\label{powerquotient}
Suppose that the complementary graph of $G$ is chordal. Then all power of $I(G)$ have linear quotients.
\end{Corollary}

To demonstrate our theory we consider the following example: let $G$ be the finite graph on the vertex set $\{1,2,3,4,5,6\}$ with edges $$\{\{1,4\},\{2,5\},\{3,6\},\{4,5\},\{4,6\},\{5,6\}\}.$$ 
The complementary graph of $G$ is chordal. 
Let $I=I(G)$ be the edge ideal of $G$, and $J$ be the toric ideal of the Rees algebra $\Rees(I)$. Then the initial ideal of $J$ with respect to the lexicographic order introduced above is generated by
\[
x_5y_1, x_4y_2, x_5y_3, x_6y_4, x_5y_5, x_4y_3, x_6y_2, x_6y_1, x_1y_4y_5, x_2y_3y_4, x_1y_3y_4, x_1y_2y_5.
\]
It follows from Corollary \ref{schappach} that $\depth S/I\geq 3$, $\depth S/I^k\geq 0$ for $k\geq 2$. Indeed, in this example  equality holds.

\medskip
\noindent
(b) Another  important class of monomial ideals with linear quotients
is the class of polymatroid ideals.
Let $I$ denote a monomial ideal of the polynomial ring $S = K[x_1, \ldots, x_n]$ generated in one degree, 
 and 
$G(I)$ its unique minimal system of monomial generators.
We say that $I$ is {\em polymatroidal} 
if the following condition is satisfied:
For monomials $u = x_1^{a_1} \cdots x_n^{a_n}$ and
$v = x_1^{b_1} \cdots x_n^{b_n}$ belonging to $G(I)$
and for each $i$ with $a_i > b_i$, one has $j$ with
$a_j < b_j$ such that $x_j u / x_i \in G(I)$.
The reason why we call such an ideal
polymatroidal is that the monomials of the ideal correspond
to the bases of a discrete polymatroid \cite{HerzogHibi}.
The polymatroidal ideal $I$ is called {\em matroidal} if
$I$ is generated by squarefree monomials.

It is known \cite[Theorem 5.2]{ConcaHerzog} that
a polymatroidal ideal has linear quotients with respect to
the reverse lexicographic order $<_{rev}$ induced by the ordering
$x_1 > x_2 > \cdots > x_n$.  
More precisely, 
if $I$ is a polymatroidal ideal and 
if $u_1, \ldots, u_s$ are the monomials belonging to $G(I)$
ordered by the reverse lexicographic order, i.e.,
$u_s <_{rev} \cdots <_{rev} u_2 <_{rev} u_1$, then
the colon ideal $(u_1, \ldots, u_{j-1}) : u_j$
is generated by a subset of $\{ x_1, \ldots, x_n \}$.

The product of polymatroidal ideals is again
polymatroidal (\cite{ConcaHerzog} and \cite{HerzogHibi}).
In particular each power of a polymatroidal ideal is 
polymatroidal.

One of the most distinguished polymatroidal ideals
is the ideal of Veronese type.  Let $S = K[x_1, \ldots, x_n]$
and fix positive integers $d$ and $e_1, \ldots, e_n$ with
$1 \leq e_1 \leq \cdots \leq e_n \leq d$.  The {\em ideal of Veronese type} 
of $S$ indexed by $d$ and $(e_1, \ldots, e_n)$ is the ideal 
$I_{(d; e_1, \ldots, e_n)}$ which is generated by 
those monomials $u = x_1^{a_1} \cdots x_n^{a_n}$ of $S$ of degree $d$ with
$a_i \leq e_i$ for each $1 \leq i \leq n$.

\begin{Theorem}
\label{depthVeronese}
Fix positive integers $d$ and $e_1, \ldots, e_n$ with
$1 \leq e_1 \leq \cdots \leq e_n \leq d$.
Let $t = d+n-1-\sum_{i=1}^{n}e_i$, and  $I = I_{(d; e_1, \ldots, e_n)}$ 
be the ideal of Veronese type of $S$
indexed by $d$ and $(e_1, \ldots, e_n)$.
Then one has $\depth S/I = t$.
\end{Theorem}

\begin{proof}
Let $u_0 = x_1^{e_1 - 1} \cdots x_{n-1}^{e_{n-1}-1} x_n^{e_n}$
and $u = x_{n-t} x_{n-t+1} \cdots x_{n-1} u_0 \in G(I)$.
For each $1 \leq i \leq n - t - 1$, one has
$x_i u / x_n \in G(I)$ with 
$u <_{rev} x_i u / x_n$.  
Let $J = (\{ w \in G(I) : u <_{rev} w \})$.
For each $1 \leq i \leq n - t - 1$, one has
$x_i u / x_n \in G(I)$ with 
$u <_{rev} x_i u / x_n$.
Hence $x_i \in J : u$
for all $1 \leq i \leq n - t - 1$.
Moreover, one has $x_j u / x_{j_0} \not\in G(I)$
for all $n - t \leq j \leq n$ and for all $j_0 \neq j$. 
Hence $x_j \not\in J : u$ for all $n - t \leq j \leq n$.
Thus $J : u = (x_1, \ldots, x_{n-t-1})$.
On the other hand, 
for each $v = x_1^{a_1} \cdots x_n^{a_n} \in G(I)$
with $m(u) = \max \{ i : a_i \neq 0 \}$,
the number of $i < m(v)$ with $a_i < e_i$ is at most 
$n - t - 1$.  
Thus the number of variables required to generate 
the colon ideal 
$(\{ w \in G(I) : v <_{rev} w \}) : v$ is at most
$n - t - 1$.
Hence $q(I) = n - t - 1$.
Thus $\depth S / I = t$.
\end{proof}

The {\em squarefree Veronese ideal} of degree $d$ in the variables
$x_{i_1}, \ldots, x_{i_t}$ is the ideal of $S$ 
which is generated by all squarefree monomials 
in $x_{i_1}, \ldots, x_{i_t}$ of degree $d$.
The squarefree Veronese ideal is matroidal and Cohen--Macaulay.

Let $2 \leq d < n$ and $I = I_{n,d}$ be
the squarefree Veronese ideal of degree $d$
in the variables $x_1, \ldots, x_n$.  Since each power $I^k$
is the ideal of Veronese type indexed by $kd$ and $(k,k,\ldots, k)$,
by using Theorem \ref{depthVeronese}, we have

\begin{Corollary}
\label{depthformula}
Let $2 \leq d < n$ and $I = I_{n,d}$ 
the squarefree Veronese ideals of degree $d$
in the variables $x_1, \ldots, x_n$.
Then
\[
\depth S / I^k = \max \{\, 0, \, n - k(n - d) - 1 \, \}.
\]
\end{Corollary}

\begin{Corollary}
\label{depthdim} 
Given nonnegative integers $d$ and $t$ with $t \leq d$
there exists a polymatroidal ideal $I \subset S$
with $\depth S / I = t$ and $\dim S / I = d$.
\end{Corollary}    

\begin{proof}
Let $I = I_{n,n-1}$ 
the squarefree Veronese ideal of degree $n-1$
in the variables $x_1, \ldots, x_{n}$.
Then $\dim S / I^k = n - 2$ and 
$\depth S / I^k = \max \{\, 0, \, n - k - 1 \, \}$.
Hence by setting $n = d + 2$ and $k = n - t - 1$,
the desired example arises.
\hspace{1.5cm}
\end{proof}

\medskip
\noindent
(c) Finally we consider a class of monomial ideals arising from finite posets. 
  Let $P$ be a finite partially ordered set ({\em poset} for short) and
write $\JJ(P)$ for the finite poset which consists of all poset
ideals of $P$, ordered by inclusion.  Here, a {\em poset ideal} of
$P$ is a subset $I \subset P$ such that if $x \in I$, $y \in P$ and
$y \leq x$, then $y \in I$.  In particular, the empty set as well as
$P$ itself is a poset ideal of $P$.  If follows that $\JJ(P)$ is a
finite distributive lattice \cite[p. 106]{StanleyEC}.  Conversely,
Birkhoff's fundamental structure theorem \cite[Theorem
3.4.1]{StanleyEC} guarantees that, for an arbitrary finite
distributive lattice $\LL$, there exists a unique poset $P$ such that
$\LL = \JJ(P)$.

Let $P = \{ p_1, \ldots, p_n \}$ be a finite poset with $|P| = n$,
and $S = K[x_1, \ldots, x_n, y_1, \ldots y_n]$ the polynomial ring in
$2n$ variables over a field $K$ with each $\deg x_i = \deg y_i = 1$. 
We associate each poset ideal $I$ of $P$ with the squarefree monomial
\[
u_I = (\prod_{p_i \in I} x_i)(\prod_{p_i \in P \setminus I} y_i)
\]
of $S$ of degree $n$.  In particular $u_P = x_1 \cdots x_n$ and
$u_\emptyset = y_1 \cdots y_n$.  We write $H_P$ for the squarefree
monomial ideal of $S$ generated by all monomials $u_I$ with $I \in
\JJ(P)$, that is, 
\[
H_P = (\{ u_I \}_{I \in \JJ(P)}).
\]

In the previous paper \cite{HerzogHibi2} it was proved that each
power $H_P^k$  has a linear resolution. 
Moreover, it is known \cite{HHZmeetsemilattice} that
$H_P$ has linear quotients.  It was expected, but unclear if all
powers of $H_P$ have linear quotients.  Fortunately, the expectation
now turns out to be true.    

\begin{Theorem}
\label{expected}
Each power $H_P^k$  has linear quotients.  
\end{Theorem} 

\begin{proof}
By virtue of \cite[p.\ 99]{HibiDistributiveLattice} each monomial
belonging to $G(H_P^k)$ possesses a unique expression of the form
$u_{I_1} u_{I_2} \cdots u_{I_k}$, where each $I_j$ is a poset ideal
of $P$, with $I_1 \subset I_2 \subset \cdots \subset I_k$.  We fix an
ordering $<$ of the monomials $u_I$, where $I$ is a poset ideal of
$P$, with the property that one has $u_I < u_J$ if $J \subset I$.  We
then introduce the lexicographic order $<_{lex}$ of the monomials
belonging to $G(H_P^k)$ induced by the ordering $<$ of the monomials
$u_I$.  We claim that $H_P^k$ has linear quotients.  More precisely,
we show that, for each monomial $w = u_{I_1} u_{I_2} \cdots u_{I_k}
\in G(H_P^k)$, the colon ideal $(\{ v \in G(H_P^k) : w <_{lex} v\}) :
w$ is generated by those variables $y_i$ for which there is $1 \leq j
\leq k$ with $p_i \in I_j$ such that $I_j \setminus \{ p_i \}$ is a
poset ideal of $P$.

First, let $y_i$ be a variable with $p_i \in I_j$ and suppose that $J
= I_j \setminus \{ p_i \}$ is a poset ideal of $P$.  One has $y_i
u_{I_j} = x_i u_J$.  Hence
\[
y_i w = x_i u_{I_1} \cdots u_{I_{j-1}} u_{J} u_{I_{j+1}} \cdots
u_{I_k}.
\]
Since each of the poset ideals $I_1, \ldots, I_{j-1}$ and $J$ is a
subset of $I_j$, it follows from 
\cite[(2.1), p.\ 98]{HibiDistributiveLattice} that the monomial $u_{I_1} \cdots
u_{I_{j-1}} u_{J}$ can be expressed uniquely in  the form $u_{I'_1}
\cdots u_{I'_{j-1}} u_{I'_j}$ such that $I'_1 \subset \cdots \subset
I'_{j-1} \subset I'_j \subset I_i$.  Moreover, one has
$u_{I_1} \cdots u_{I_{j-1}} u_{J} <_{lex} 
u_{I'_1} \cdots u_{I'_{j-1}} u_{I'_j}$.
Thus $w <_{lex} 
u_{I'_1} \cdots u_{I'_{j-1}} u_{I'_j} u_{I_{j+1}} \cdots u_{I_k}$.
Hence $y_i$ belongs to the colon ideal 
$(\{ v \in G(H_P^k) : w <_{lex} v\}) : w$.

Second, let $\delta$ be a monomial belonging to the colon  ideal
$$(\{ v \in G(H_P^k) : w <_{lex} v\}) : w.$$  Thus one has $\delta w =
\mu v$ for monomials $\mu$ and $v$ with $w <_{lex} v$.  Say, 
$v = u_{I'_1} \cdots u_{I'_k}$ with $I'_1 \subset \cdots \subset
I'_k$.
What we must prove is that the monomial $\delta$ is divided by  a
variable $y_i$ for which there is $1 \leq j \leq k$ such that $I_j
\setminus \{ p_i \}$ is a poset ideal of $P$.  Since $w <_{lex} v$,
it follows that there is $j_0$ for which
$I_{j_0} < I'_{j_0}$.  In particular $I_{j_0} \not\subset I'_{j_0}$. 
Thus there  is a maximal element $p_{i_0}$ of $I_{j_0}$ with $p_{i_0}
\not\in I'_{j_0}$.  Then $p_{i_0}$ belongs to each of the poset
ideals $I_{j_0}, I_{j_0 + 1}, \ldots, I_k$ and belongs to {\em none}
of the poset ideals $I'_1, \ldots, I'_{j_0}$. Hence the power of
$y_{i_0}$ in the monomial $v$ is at least $j_0$, but that in $w$ is
at most $j_0 - 1$.  Hence $y_0$ must divide $\delta$.  Since $p_{i_0}$ is
a maximal element of $I_{j_0}$, the subset $I_{j_0} \setminus \{
p_{i_0} \}$ of $P$ is a poset ideal of $P$, as desired. 
\end{proof}

By using Theorem \ref{expected} we can now compute $\depth S/H_P^k$
in terms of the combinatorics on $P$.  Recall that an {\em antichain}
of $P$ is a subset $A \subset P$ any two of whose elements are
incomparable in $P$.  Given an antichain $A$ of $P$, we write $\langle
A \rangle$ for the poset ideal of $P$ generated by $A$, which consists of those elements $p
\in P$ such that there is $a \in A$ with $p \leq a$.  For each $k =
1, 2, \ldots$, we write $\delta(P;k)$ for the largest integer $N$ for
which there is a sequence $(A_1, A_2, \ldots, A_r)$ of antichains of $P$ with $r \leq k$  such
that
\begin{enumerate}
\item[(i)]
$A_i \cap A_j = \emptyset$ if $i \neq j$;
\item[(ii)]
$\langle A_1 \rangle \subset 
\langle A_2 \rangle \subset 
\cdots 
\subset \langle A_r \rangle$;  
\item[(iii)] $N = |A_1| + |A_2| +, \cdots + |A_r|.$
\end{enumerate} 
We call such a sequence of antichains  a {\em
$k$-acceptable sequence}. 
 
It follows from the definition that  $\delta(P;1)$ is the maximal cardinality of
antichains of $P$ and $\delta(P;1)< \delta(P;2) < \cdots <
\delta(P;\rank(P) + 1)$.  Moreover, $\delta(P;k) = n$ for all $k \geq
\rank(P) + 1$.  Here $\rank(P)$ is the {\em rank} \cite[p.\  99]{StanleyEC} of $P$.  Thus $\rank(P) + 1$ is the maximal
cardinality of chains (totally ordered sets) contained in $P$.

\begin{Corollary}
\label{posetdepth} 
Let $P$ be an arbitrary finite poset with $|P| = n$.  Then  
\[
\depth S/H_P^k = 2n - \delta(P;k) - 1
\]
for all $k \geq 1$.
\end{Corollary}

\begin{proof}
We work with the same notation as in the proof of Theorem
\ref{expected}.  Recall that, for a monomial $w = u_{I_1} u_{I_2}
\cdots u_{I_k} \in G(H_P^k)$, the colon ideal $(\{ v \in G(H_P^k) : w
<_{lex} v\}) : w$ is generated by those variables $y_i$ for which
there is $1 \leq j \leq k$ with $p_i \in I_j$ such that $I_j
\setminus \{ p_i \}$ is a poset ideal of $P$.  Note that $I_j
\setminus \{ p_i \}$ is a poset ideal of $P$ if and only if $p_i$ is
a maximal element of $I_j$.  Let $B_j$ denote the set of maximal
elements of $I_j$.  Then
the number of variables required to generate the colon ideal $(\{ v
\in G(H_P^k) : w <_{lex} v\}) : w$ is $|\Union_{j=1}^{k} B_j|$.  Let
$Q_w = \Union_{j=1}^{k} B_j$.  One has $r = \rank(Q_w) + 1 \leq k$. 
We then define a sequence $A_1, A_2, \ldots, A_{r}$ of subset of
$B_w$ as follows:  $A_1$ is the set of minimal elements of $Q_w$ and,
for $2 \leq j \leq r$, $A_j$ is the set of minimal element of $Q_w
\setminus (A_1 \cup \cdots \cup A_{j-1})$.  Then $(A_1, \ldots, A_r)$ 
is $k$-acceptable with 
$|Q_w| = \sum_{j=1}^{r} |A_j|$.  Hence $|Q_w| \leq \delta(P;k)$.  

On the other hand, there is a $k$-acceptable sequence $(A_1, A_2,
\ldots, A_r)$ with $\delta(P;k) = \sum_{j=1}^{r} |A_j|$.  Let 
$w = u_{\emptyset}^{k - r} u_{\langle A_1 \rangle} \cdots u_{\langle
A_r \rangle} \in G(H_P^k)$.  Then the number of variables required to
generate the colon ideal $(\{ v \in G(H_P^k) : w <_{lex} v\}) : w$ is
$\delta(P;k)$.

Consequently, one has $q(H_P^k) = \delta(P;k)$.  Thus $\depth S/H_P^k
= 2n - \delta(P;k) - 1$, as required.
\end{proof}

Since $\{ x_i, y_i \}$ is a vertex cover of $H_P$ for each $1 \leq i
\leq n$, it follows that $\dim S/H_P = 2n - 2$.  Hence $H_P$ is
Cohen--Macaulay if and only if $\delta(P;1) = 1$.  In other words,
$H_P$ is Cohen--Macaulay if and only if $P$ is a chain.

\begin{Corollary}
\label{posetdepthlimit} 
Let $P$ be an arbitrary finite poset with $|P| = n$.  Then 
\begin{enumerate}
\item[(i)]
$\depth S/H_P > \depth S/H_P^2 > \cdots > \depth S/H_P^{\rank(P)} 
> \depth S/H_P^{\rank(P) + 1}$;
\item[(ii)] 
$\depth S/H_P^{k} = n - 1$ for all $k > \rank(P)$;
\item[(iii)]  
$\lim_{k \to \infty} \depth S/H_P^k = n - 1$.
\end{enumerate}
\end{Corollary}

\begin{Corollary}
\label{directsum}
Given an integer $n > 0$ and given a finite sequence $(a_1,
a_2, \ldots, a_r)$ of positive integers with $a_1 \geq a_2 \geq
\cdots \geq a_r$ and with $a_1 + \cdots + a_r = n$, there exists a
squarefree monomial ideal $I \subset S = K[x_1, \ldots, x_n, y_1,
\ldots, y_n]$ such that
\begin{enumerate}
\item[(i)]
$\depth S/I^k = 2n - (a_1 + \cdots + a_k) - 1$, \, \, \, \, \, $k = 1,
2, \ldots, r - 1$;
\item[(ii)] 
$\depth S/I^k = n - 1$ for all $k \geq r$;
\item[(iii)]
$\lim_{k \to \infty} \depth S/I^k = n - 1$.
\end{enumerate}
\end{Corollary}

\begin{proof}
Let $A(a_i)$ denote the antichain with $|A(a_i)| = a_i$ and $P$ the
ordinal sum \cite[p.\ 100]{StanleyEC} of the antichains $A(a_1),
A(a_2), \ldots, A(a_r)$.
Thus $\rank(P) = r - 1$.  
Since $a_1 \geq a_2 \geq \cdots \geq a_r$ and $a_1 + \cdots + a_r =
n$, it follows that $\delta(P;k) = a_1 + a_2 + \cdots + a_k$ if $1
\leq k \leq r - 1$
and that $\delta(P;k) = n$ for all $k \geq r$.     
\end{proof}

In general, given a function $f : \NN  \to \NN$,  we
introduce the function $\Delta f$ by setting $(\Delta f)(k) = f(k) -
f(k+1)$ for all $k \in \NN$. 

\begin{Corollary}
\label{decreasing}
Given a decreasing function $f : \NN \to \NN$ with 
$$f(0) = 2 \lim_{k \to \infty} f(k) + 1$$ 
for which $\Delta f$ is decreasing, there exists a monomial ideal $I
\subset S$ such that $\depth S / I^k\newline = f(k)$ for all $k \geq 1$. 
\end{Corollary}

\begin{proof}
Let $\lim_{k \to \infty} f(k) = n - 1$ and $f(0) = 2n -1$.  
Let $a_k = (\Delta f)(k - 1)$ for all $k \geq 1$.
Thus $f(k) = 2n - (a_1 + \cdots a_k) - 1$ for all $k \geq 1$.
Since $f$ is decreasing, one has $a_k \geq 0$ for all $k$.
Since $\Delta f$ is decreasing, one has $a_1 \geq a_2 \geq \cdots$.
Let $r \geq 1$ denote the smallest integer for which 
$a_1 + a_2 + \cdots + a_r = n$.
Thus $a_i > 0$ for $1 \leq i \leq r$ and $a_i = 0$ for all $i > r$.
It then follows from Corollary \ref{directsum} that
there exists a monomial ideal $I \subset S$ for which
$\depth S / I^k = f(k)$ for all $k \geq 1$.
\end{proof}

\section{A class of ideals whose depth function $\depth S/I^k$ is increasing}

Note that if $I$ is a squarefree  monomial ideal, then $\depth S/I^k\leq \depth S/I$ for all $k$, see for example \cite{HerzogTakayamaTerai}. This suggests the following question: Is it true that $\depth S/I^k$ is a decreasing function of $k$, if $I$ is a squarefree monomial ideal?
As we shall see now, for a general monomial ideal the function $\depth S/I^k$ may also be increasing. In fact, we even show 

\begin{Theorem}
\label{bad}
Given a bounded increasing function $f\: \NN\setminus\{0\}\to \NN$. There exists a monomial ideal $I$ such that $\depth S/I^k=f(k)$ for all $k$.
\end{Theorem}

\begin{proof}  Let $\lim_{k\to\infty}f(k)=n$, and suppose that $f(k)=n$ for $k\geq d-1$. We set 
\begin{eqnarray}
\label{definition}
c_{d-k}=n-f(k)\quad \text{for}\quad k=1,\ldots, d-2.
\end{eqnarray}
Let $K$ be field, and $S=K[x_1,x_2,y_1,\ldots, y_n]$ be the polynomial ring in $n+2$ variables over $K$.  
We define  $I\subset S$ to be the ideal generated by the set of monomials
\[
\{x_1^{d+1},x_1^dx_2,x_1x_2^d,x_2^{d+1}\}\union\Union_{k=2}^{d-1}\{x_1^{d-1}x_2^ky_1,\ldots,  x_1^{d-1}x_2^ky_{c_k}\}.
\]
Note that this set of monomials is in general not a minimal set of generators of $I$. We claim that 
\[
\depth S/I^k=f(k)\quad \text{for all $k$.}
\] 
For $k=1,\ldots, d-2$, let $J_{(k)}\subset S_r=K[x_1,x_2,y_1,\ldots, y_{c_{d-k}}]$ be the ideal generated by the set of monomials
\[
\{x_1^{d+1},x_1^dx_2,x_1x_2^d,x_2^{d+1}\}\union\Union_{r=2}^{d-k}\{x_1^{d-1}x_2^ry_1,\ldots,  x_1^{d-1}x_2^ry_{c_r}\},
\]
and set $J=J_{(d-1)}= (x_1^{d+1},x_1^dx_2,x_1x_2^d,x_2^{d+1})$. We will show:
\begin{enumerate}
\item[(i)] $J_{(k)}^kS=I^k$ for $k=1,\ldots, d-1$, 
\item[(ii)] $x_1^{kd-1}x_2^{d-1}\not\in J_{(k)}^k$ for $k=1,\ldots, d-2$, and
\item[(iii)] $x_1^{kd-1}x_2^{d-1}(x_1,x_2,y_1,\ldots, y_{d-k})\in J_{(k)}^k$. 
\end{enumerate}
Assuming (i), (ii) and (iii), the assertion follows. Indeed, if we set $c_1=0$, then (i) implies
\begin{eqnarray}
\label{deptheq}
\depth S/I^k=\depth S/J_{(k)}^kS=\depth S_k/J_{(k)}^k+ (n-c_{d-k}),
\end{eqnarray}
for $k=1,\ldots, d-1$, and (ii) and (iii) imply that $\depth S_k/J_{(k)}^k=0$. Thus (\ref{definition}) and (\ref{deptheq}) yield the desired result.

Before proving (i), (ii) and (iii) we notice that $J^k$ is generated in degree $k(d+1)$, and that for any $r\geq 
k(d+1)$ one has
\begin{eqnarray}
\label{elements}
(J^k)_{\geq r}&=&[(x_1^d,x_2^d)^k(x_1,x_2)^k](x_1,x_2)^{r-k(d+1)}=(x_1^d,x_2^d)^k(x_1,x_2)^{r-kd}\\
&=& (\{x_1^{id+s}x_2^{r-(id+s)}\}_{i=0,\ldots,k,\atop s=0,\ldots,r-kd}).\nonumber
\end{eqnarray}
Proof of (i): The desired equality follows once we can show for all $t=1,\ldots, k$ the ideal $J^{k-t}$ multiplied with a product of $t$ elements from the set $$\Union_{k=2}^{d-1}\{x_1^{d-1}x_2^ky_1,\ldots,  x_1^{d-1}x_2^ky_{c_k}\},$$ with at least one factor of the form $x_1^{d-1}y_2^ry_i$ with $r\geq d-k+1$, belongs to $J^k$. This will be the case if $J^{k-t}(x_1^{d-1}x_2^{r_1})\cdots (x_1^{d-1}x_2^{r_t})\subset J^k$ for all $t=1,\ldots,k$ and all $r_i$ with $2\leq r_1\leq r_2\leq \cdots \leq r_t$ with at least one $r_i\geq d-k+1$. For this it suffices to consider the most critical case, namely that $r_1=r_2=\cdots,r_{t-1}=2$ and $r_t=d-k+1$. Thus we have to show that $J^{k-t}(x_1^{d-1}x_2^2)^{t-1}(x_1^{d-1}x_2^{d-k+1})\subset J^k$. By (\ref{elements}) it amounts therefore to show that 
\[
u=x_1^{id+s}x_2^{(k-t)(d+1)-(id+s)}x_1^{t(d-1)}x_2^{2t+d-k-1}\in (J^k)_r,
\]
where $r=(k-t)(d+1)+(t-1)(d-1)+2(t-1)+d-1+d-k+1=kd+d-1$ is the degree of the monomial $u$, and where $0\leq i\leq k-t$ and $0\leq s\leq k-t$. Again using (\ref{elements}) we see that $u\in (J^k)_r$ if and only if 
\[
(i+t)d+(d-t+s)\in \{jd+a\: 0\leq j\leq k,\; 0\leq a\leq d-1\}.
\]
Since $0\leq i\leq k-t$ it follows that $1\leq i+t\leq k$ is in the allowed range. We also have
\[
2\leq d-k\leq d-t+s\leq 2d-1.
\]
If $i+t=k$, then $t=k$ and $s=0$, so that  $d-t+s=d-k\leq d-1$. On the other hand, if $i+t<k$, then 
$(i+t)d+(d-t+s)= (i+t+1)d+(-t+s)$ has the desired form.

\medskip
\noindent
Proof of (ii): It suffices to show that $x_1^{kd-1}x_2^{d-1}\not\in J^k$, because the ideals $J^k$ and $J_{(k)}^k$ coincide modulo $y_1,\ldots, y_n$. 

Suppose $x_1^{kd-1}x_2^{d-1}\in J^k$, then $x_1^{kd-1}x_2^{d-1}\in (J^k)_{kd+(d-2)}$. It follows from (\ref{elements}), that $$kd-1\in \{id+s\: i=0,\ldots,k,\; \text{and  $s=0,\ldots,d-2$}\}.$$
Hence we must have $kd-1=id+s$ for some $0\leq i<k$. This yields $(k-i)d=s+1\leq d-1$, a contradiction.  

\medskip
\noindent
Proof of (iii): The element  $x_1(x_1^{kd-1}x_2^{d-1})=x_1^{kd}x_2^{d-1}$ belongs  to $(J^k)_{kd+(d-1)}\subset J_{(k)}^k$. Also, by (\ref{elements}),  the element $x_2(x_1^{kd-1}x_2^{d-1})=x_1^{kd-1}x_2^{d}$ belongs to $(J^k)_{kd+(d-1)}$, since $kd-1=(k-1)d+(d-1)$. Finally, we note that for $i=1, \ldots,c_{d-k}$ one has
\[
y_i(x_1^{kd-1}x_2^{d-1})= (x_1^{(k-1)d}x_2^{k-1})(x_1^{d-1}x_2^{d-k}y_i).
\]
By (\ref{elements}), the first factor belongs to $(J^{k-1})_{(k-1)(d+1)}$, and the second factor belongs to 
$J_{(k)}$. Thus $y_i(x_1^{kd-1}x_2^{d-1})=J_{(k)}^k$, as desired.
\end{proof}

All examples we have considered so far had the property that the function $\depth S/I^k$ is monotonic. We conclude this paper with an example that shows that this depth function can be more general. We consider the ideal
\[
I=(a^6,a^5b,ab^5,b^6,a^4b^4c,â^4b^4d,a^4e^2f^3,b^4e^3f^2)\
\]
in $S=K[a,b,c,d,e,f]$. Then $\depth S/I=0$, $\depth S/I^2=1$, $\depth S/I^3=0$,  $\depth S/I^4=2$ and $\depth S/I^5=2$.

In view of the examples considered in this paper we are   tempted to conjecture that the function $\depth S/I^k$ can be {\em any} convergent nonnegative integer valued function.

\end{document}